\documentclass{article}
\usepackage{amsmath,amsfonts,amssymb}
\usepackage[latin2]{inputenc}                                          

\def\softd{{\leavevmode\setbox1=\hbox{d}%
\hbox to 1.05\wd1{d\kern-0.4ex{\char039}\hss}}}
\def\softt{{\leavevmode\setbox1=\hbox{t}
\hbox to \wd1{t\kern-0.6ex{\char039}\hss}}}
\def\softl{l\kern-0.45ex\raise0.1ex\hbox{''}\kern-0.10ex}
\def\softL{L\kern-0.8ex\raise0.1ex\hbox{''}\kern0.1ex}

\newtheorem{theorem}{Theorem}[section]
\newtheorem{definition}[theorem]{Definition}

\newtheorem{corollary}[theorem]{Corollary}

\newtheorem{example}[theorem]{Example}
\newtheorem{lemma}[theorem]{Lemma}

\newtheorem{remark}[theorem]{Remark}

\newcommand{\la}{\mathbb \lambda}

\newcommand{\mY}{\mathbb Y}

\newcommand{\mV}{\mathbb V}
\newcommand{\mH}{\mathbb H}
\newcommand{\mC}{\mathbb C}

\newcommand{\lD}{\cal D}

\newcommand{\mZ}{\mathbb Z}

\newcommand{\mR}{\mathbb R}
\newcommand{\mW}{\mathbb W}
\newcommand{\mN}{\mathbb N}
\newcommand{\mE}{\mathbb E}

\newcommand{\be}{\begin{eqnarray}}
\newcommand{\ee}{\end{eqnarray}}
\newcommand{\bd}{\begin{definition}}
\newcommand{\ed}{\end{definition}}
\newcommand{\br}{\begin{remark}}
\newcommand{\er}{\end{remark}}

\newcommand{\gog}{{\mathfrak g}}
\newcommand{\gol}{{\mathfrak l}}

\newcommand{\gon}{{\mathfrak n}}
\newcommand{\gop}{{\mathfrak p}}

\newcommand{\bt}{\begin{tabular}}
\newcommand{\et}{\end{tabular}}

\newcommand{\bl}{\begin{lemma}}
\newcommand{\el}{\end{lemma}}

\newcommand{\bp}{\begin{picture}}
\newcommand{\ep}{\end{picture}}
\newcommand{\bi}{\begin{itemize}}
\newcommand{\ei}{\end{itemize}}
\newcommand{\bq}{\begin{quotation}}
\newcommand{\eq}{\end{quotation}}

\newcommand{\Pol}{\operatorname{Pol}}

\newcommand{\Ind}{\operatorname{Ind}}
\newcommand{\Hom}{\operatorname{Hom}}
\newcommand{\ch}{{\,\vee}}
\newcommand{\tpi}{\tilde{\pi}}
\newcommand{\Id}{\operatorname{Id}}

\input epsf
\input xy
\xyoption{all}

\begin{document}
\baselineskip13pt

\title{Rankin-Cohen brackets for orthogonal Lie algebras and 
bilinear conformally equivariant differential operators
}

\author{Petr Somberg}

\date {}

\maketitle

\center{\it Dedicated to Professor Toshi Kobayashi on the occasion of his birthday,}
\center{\it with admiration.}
\vspace{1cm}

\abstract 
Based on the Lie theoretical methods of algebraic Fourier transformation,
we classify in the case of generic values of inducing parameters the scalar singular 
vectors corresponding to the diagonal branching rules for scalar generalized 
Verma modules in the case of orthogonal Lie algebra and its conformal parabolic 
subalgebra with commutative nilradical, thereby realizing the diagonal branching 
rules in an explicit way. The complicated combinatorial structure of singular vectors is conveniently 
determined in terms of recursion relations for the generalized hypergeometric function 
${}_3F_2$. 
As a geometrical application, we classify bilinear conformally equivariant differential 
operators acting on homogeneous line bundles on the flag manifold given by conformal sphere $S^n$.    

\hspace{0.2cm}

 

\endabstract

\section{Introduction and Motivation}

The subject of our article has its motivation in the Lie theory 
for finite dimensional simple Lie algebras applied to the problem of 
branching rules and composition structure of generalized Verma modules, 
and dually in the geometrical situation related to the construction of equivariant bilinear 
differential operators or Rankin-Cohen-like brackets associated to orthogonal groups.    

The classical Rankin-Cohen brackets
realized by holomorphic $SL(2,\mR)$-equivariant bilinear differential operators 
on the upper half plane $\mH$ are devised, originally in a number theoretic context,
to produce from a given pair of modular forms another modular form.
They turn out to be intertwining operators responsible for the ring structure on 
$SL(2,\mR)$ holomorphic discrete series representations, and can 
be analytically continued to the full range of inducing characters. Consequently,
such operators were constructed by different techniques in several specific situations 
of interest related to Jacobi forms, Siegel modular forms, holomorphic discrete 
series of causal symmetric spaces of Cayley type, real symmetric pairs of split 
rank one, etc., cf. \cite{grad}, \cite{ey}, \cite{dp}, \cite{kope}.

The main result of the present article is the classification of Rankin-Cohen 
conformally equivariant bilinear differential (bidifferential) operators acting 
on sections of homogeneous
line bundles over the conformal sphere $S^n$. These operators can be regarded as 
projectors onto irreducible summands in the decomposition of the tensor product 
of two particular representations of $SO_0(n+1,1,\mR)$. 
Our approach to this geometrical problem for finding equivariant bilinear differential operators
is based on its conversion to a Lie algebraic problem of the characterization of homomorphisms 
of generalized Verma modules. Based on the techniques of algebraic Fourier transform for generalized
Verma modules (F-method) recently initiated in \cite{koss}, \cite{kope1}, \cite{kope}, we solve this problem in the case of
orthogonal Lie algebras and their conformal (commutative) parabolic subalgebras. The finer 
questions related to the description of composition series for special values of infinitesimal characters
or the complete classification of all solutions for the system of PDEs via F-method are postponed 
to a future research. 

An abstract description of the underlying multiplicity-free branching problem 
in question is known, cf. \cite{ko}. Our intention is the computation of precise positions of submodules 
for the branching subgroup in the whole representation space, whose explicit knowledge is 
required in many applications like geometric analysis on manifolds, number theory, etc.

There are various approaches to the questions discussed in our short article. For example,
there is an analytic approach consisting of meromorphic continuation of invariant 
distributions given by a multilinear form on the principal series representations. A class of 
$SO_0(n+1,1,\mR)$ (i.e., conformally)-covariant linear and bilinear differential
operators was realized in residues of meromorphically continued invariant 
trilinear form on principal series representations induced from characters, 
see e.g., \cite{bc}.
In the geometrical context of generalization towards curved manifolds with parabolic 
structure, a classification of first order equivariant bilinear differential operators for parabolic 
subalgebras with commutative nilradicals (AHS structures) was completed in \cite{kroes}.
As for a geometrical application of conformally covariant bilinear differential operators,
see \cite{me}.
We remark that the appearance of a higher hypergeometric functions in our main 
Theorem \ref{bilconfdensclass} parallels the results in \cite{kope} on Rankin-Cohen
bracket for real split rank one symmetric spaces, e.g. the case of $(U(n,1)\times U(n,1),U(n,1))$ 
corresponding to the diagonal embedding of complex projective space
$\mC{\mathbb P}^n\hookrightarrow \mC{\mathbb P}^n\times\mC{\mathbb P}^n$ where the  
differential operators in concern are given by substituting
differential operators into appropriate hypergeometric function $_{2}F_1$
(Jacobi and Gegenbauer polynomials.)
 
The structure of our article goes as follows. In Section $2$, we reformulate the question of 
existence for equivariant bilinear differential operators (the Rankin-Cohen brackets) in terms 
of purely abstract Lie theoretic classification scheme for diagonal branching rules of 
generalized Verma modules, and apply it to the case of real orthogonal Lie algebras $so(n+1,1,\mR)$ and 
their conformal parabolic Lie subalgebras $\gop$. We review some necessary technical background,
and discuss the abstract (or, qualitative) part of the branching problem taking its value in 
the Grothendieck group $K({\fam2 O}^\gop)$ of the Bernstein-Gelfand-Gelfand parabolic category 
${\fam2 O}^\gop$. We also review the procedure of algebraic Fourier transformation applied to 
generalized Verma modules, which will be used in Section $3$.   
The main theme of our article is the quantitative part of the branching problem, see Section $3$,
which consists of the construction of a class of (scalar valued) singular vectors. Our approach to analyze singular vectors
is based on the procedure of F-method applied to generalized Verma modules, 
where the action of positive nilradical of $so(n+1,1,\mR)$ reduces to a four term functional equation 
for the coefficients of singular vectors. To solve it is technically the most 
difficult part, including both analytic and combinatorial aspects of 
generalized hypergeometric functions ${}_3F_2$. 
The results of Sections $2,3$ are applied in Section $4$ to produce classification
and explicit formulas of bilinear conformally covariant differential operators 
corresponding to (scalar valued) singular vectors.

There are several ways allowing to produce the results equivalent
to ours. For example, \cite{OvRe} is based on the construction of
transvectants, \cite{Cl1} and \cite{Cl2} develop the (meromorphic
continuation of) conformally invariant trilinear forms, or even 
the construction of (explicit formulas are of low order only) curved 
analogues on a Riemannian manifold via conformally invariant ambient 
metric construction in \cite{CaLY}. It follows from our construction 
via F-method the completeness of the constructed set, a property which
is not automatic in the other approaches. In addition, our results do agree
with those appearing in the above mentioned references. Moreover, 
an additional effort can be used to construct the lifts to homomorphisms 
of semi-holonomic Verma modules in conformal geometry, which then give 
the curved version 
of our construction parallel to the results in \cite{CaLY}. 

Throughout the article we denote by $\mN$ the set of natural numbers including zero.

\section{F-method and diagonal branching problem for generalized Verma modules}

In the present section we briefly review basic notations and results initiated and 
developed in \cite{kob}, \cite{kope}, \cite{koss}, \cite{kope1}, \cite{ko}, allowing in an explicit 
way to realize the diagonal branching problem for 
real orthogonal Lie algebras and their conformal parabolic subalgebras.

We denote by ${G}_{\mathbb R}$ a connected real reductive Lie group
 with real Lie algebra $\gog_{\mathbb R}$, $P_{\mathbb R}\subset G_{\mathbb R}$ 
 a parabolic subgroup and $\gop_{\mathbb R}$ its Lie algebra, 
 $\gop_{\mathbb R}=\gol_{\mathbb R}\oplus\gon_{\mathbb R}$
the Levi decomposition of $\gop_{\mathbb R}$ and ${\gon_{\mathbb R}}_-$ its opposite nilradical,
$\gog_{\mathbb R}={\gon_{\mathbb R}}_-\oplus\gop_{\mathbb R}$. 
Given a complex finite dimensional $P_{\mathbb R}$-module $\mV$, we consider the induced 
representation $(\pi,{\Ind}_{P_{\mathbb R}}^{G_{\mathbb R}}(\mV))$ of $G_{\mathbb R}$ on the 
space of smooth complex valued sections of the homogeneous 
vector bundle $G_{\mathbb R} \times_{P_{\mathbb R}} \mV \to G_{\mathbb R}/ P_{\mathbb R}$, 
where ${\Ind}_{P_{\mathbb R}}^{G_{\mathbb R}}(\mV):=C^\infty(G_{\mathbb R},\mV)^{P_{\mathbb R}}$ with
\begin{eqnarray}
C^\infty(G_{\mathbb R},\mV)^{P_{\mathbb R}}=\{f\in C^\infty(G_{\mathbb R},\mV)|
f(g\cdot p)=p^{-1}\cdot f(g),\, g\in G_{\mathbb R}, p\in P_{\mathbb R} \}.
\end{eqnarray}
We denote the complexification of a given real Lie algebra or group by omitting 
its subscript ${\mathbb R}$, so for example $\gog=\gog_{\mathbb R}\otimes_\mR\mC$, etc. Then 
${\fam2 U}(\gog)$ denotes the universal enveloping algebra of ${\gog}$.
For $\mV^\ch$ the contragradient (dual) complex representation to $\mV$,
the generalized Verma module ${\fam2 M}^{\gog}_{\gop}(\mV^\ch)$ 
is defined by
\begin{eqnarray}
{\fam2 M}^{\gog}_{\gop}(\mV^\ch):= 
{\fam2 U}(\gog)\otimes_{{\fam2 U}(\gop)}\mV^\ch .  
\end{eqnarray}
 There is a $({\mathfrak {g}}, {P}_{\mathbb R})$-invariant natural pairing 
\begin{equation}\label{distrduality}
\Ind^{G_{\mathbb R}}_{P_{\mathbb R}}(\mV)\times{\fam2 M}^{\gog}_{\gop}(\mV^\ch) \longrightarrow \mC,
\end{equation}
where the space of $G_{\mathbb R}$-equivariant differential 
operators $\Ind_{P_{\mathbb R}}^{G_{\mathbb R}}(\mV)\to\Ind_{P_{\mathbb R}}^{G_{\mathbb R}}(\mV')$ 
is bijective to the space of $(\gog, P_{\mathbb R})$-homomorphisms 
${\fam2 M}^{\gog}_{\gop}(\mV'^\ch)\to{\fam2 M}^{\gog}_{\gop}(\mV^\ch)$. 
Recall that homomorphisms of generalized Verma modules are determined by their singular vectors.

A generalization of the previous framework is based on two compatible pairs of real Lie 
groups $(G_{\mathbb R},P_{\mathbb R})$ and 
$( G_{\mathbb R}',P_{\mathbb R}')$, where $G_{\mathbb R}'\subset G_{\mathbb R}$ is a real reductive subgroup 
of $G_{\mathbb R}$ and 
$P_{\mathbb R}'=P_{\mathbb R} \cap G_{\mathbb R}'$ is compatible parabolic subgroup of $G_{\mathbb R}'$. The Lie algebras
of $G_{\mathbb R}',P_{\mathbb R}'$ are denoted by ${\mathfrak {g}}_{\mathbb R}',{\mathfrak {p}}_{\mathbb R}'$, 
${\mathfrak {n}}_{\mathbb R}':={\mathfrak {n}}_{\mathbb R} \cap {\mathfrak {g}}_{\mathbb R}'$ is the nilradical 
 of ${\mathfrak {p}}_{\mathbb R}'$, 
 and $L_{\mathbb R}'=L_{\mathbb R} \cap G_{\mathbb R}'$
 is the Levi subgroup of $P_{\mathbb R}'$. As in the previous paragraph, omitting
the subscript ${\mathbb R}$ denotes the complexification of real Lie algebra or group. 
Therefore, an irreducible $L_{\mathbb R}'$-submodule $\mW^{\ch}$ of 
\begin{eqnarray}
{\fam2 M}_{\gop}^{\gog}(\mV^\ch)^{\gon'}
 :=
 \{v\in {\fam2 M}^{\gog}_{\gop}(\mV^\ch)|\, Z\cdot v=0
\textrm{  for all } Z\in\gon'\}
\end{eqnarray}
with the standard left action of $Z$ on the generalized Verma module, 
gives a ${\fam2 U}(\gog')$-homomorphism
${\fam2 M}_{\gop'}^{\gog'}(\mW^\ch)\to{\fam2 M}_{\gop}^{\gog}(\mV^\ch)$. 

The whole procedure of the F-method to find explicit singular 
vectors may be divided into the following three main steps:

\vskip 1mm
 \noindent
{\bf Step 1.} Computation of the infinitesimal action $d\pi(X)$ for $X\in \gon_{\mathbb R}$
 on a chosen principal series representation $(\pi, \Ind^{G_{\mathbb R}}_{P_\mathbb R}(\mV))$ of 
 $G_{\mathbb R}$, realized in the non-compact picture.
 The induced representation $d\pi$ defines, by its restriction, the representation of $\gog$ on 
 the space ${C}^\infty({N_{\mathbb R}}_-,\mV).$

\vskip 1mm 
 \noindent
{\bf Step 2. }
Computation of the dual infinitesimal action $d\pi^\ch(X)$ for $X\in \gon_{\mathbb{R}}$ on the dual space 
${\lD}'_{[o]}({N_{\mathbb R}}_-, \mV^\ch)$
of distributions
on ${N_{\mathbb R}}_-$ with values in $\mV^\ch$ and supported in the unit $P_{\mathbb R}$-coset $[o].$
Here we recall that we realize ${N_{\mathbb R}}_-$ as the open Bruhat cell ${N_{\mathbb R}}_-P_{\mathbb R}$ of
the flag manifold $G_{\mathbb R}/P_{\mathbb R}$.
The Lie algebra $\gog$
acts on the space of vector valued distributions by the dual action $d\pi^\ch:$
$$
d\pi^\ch(X)(T)(f)=-T(d\pi(X)(f)), \quad X\in\gog,\, f\in{C}^\infty({N_{\mathbb R}}_-,\mV).
$$

This space is isomorphic with $M^{\gog}_{\gop}(\mV^\ch)$ as $\gog$-modules 
by linear map
$$
\phi:{\fam2 U}(\gog)\otimes_{{\fam2 U}(\gop)}\mV^\ch\to {\lD}'_{[o]}({N_{\mathbb R}}_-,\mV^\ch)
$$
determined by
$$
\phi(u\otimes v^\ch):f\mapsto \langle v^\ch, (d\pi(u^o)f)(o) \rangle
$$
with the map $u\to u^o$ given by the antiautomorphism of ${\fam2 U}(\gog)$ 
acting as $X\mapsto -X$ on $\gog.$

\vskip 1mm
\noindent
{\bf Step 3.}
Let us identify ${\gon_{\mathbb R}}_-$ with ${N_{\mathbb R}}_-$ by the exponential 
map. Based on the convention introduced in \cite{koss}, the Fourier transform 
gives an isomorphism
\begin{eqnarray}\label{eqn:FT}
{\mathcal F}: {\lD}'_{[0]}({\gon_{\mathbb R}}_-)
 \stackrel{\sim}{\longrightarrow} \Pol[\mathfrak {n}],\quad T\mapsto {\mathcal F}(T)
\end{eqnarray}
defined by
$$
{\mathcal F}(T)(\xi)=T_x(e^{{\mathbf i}\langle x,\xi\rangle}),\,
\mbox{for}\,\,\, x\in{\gon_{\mathbb R}}_-,\,\xi\in\gon,
$$
where $\langle x,\xi\rangle$ is given by the Killing form on $\gog$ and 
${\mathbf i}\in{\mathbb C}$ denotes the complex unit. 
The isomorphism (\ref{eqn:FT}) can be extended to distributions with values 
 in $\mV^\ch$ by
 $$
 {\mathcal F}\otimes \Id_{\mV^\ch}:{\lD}'_{[0]}({{\mathfrak n}_{\mathbb R}}_-,\mV^\ch)
 \to \Pol[\mathfrak {n}]\otimes \mV^\ch.
 $$
 The Fourier transform $ {\mathcal F}\otimes \Id_{\mV^\ch}$ can be then used 
to define the action $d\tilde{\pi}^\ch$
of $\gog$ on the space $\Pol[\mathfrak {n}]\otimes \mV^\ch$ by
$$
d\tilde{\pi}^\ch = ({\mathcal F}\otimes \Id_{\mV^\ch})\circ d{\pi}^\ch
\circ ({\mathcal F}^{-1}\otimes \Id_{\mV^\ch}),
$$
where the
elements of $\gon$ act  by differential operators of second order as ${\gon}$ 
is commutative. We set
$\varphi:=\phi^{-1}\circ ({\mathcal F}^{-1}\otimes \mathrm{Id}_{\mV^\ch})$.
Then $\varphi$ gives a bijection 
\begin{eqnarray}\label{eqn:phiPM}
\varphi: \Pol[\mathfrak {n}]\otimes \mV^\ch
\stackrel{\sim}{\longrightarrow} 
{\fam2 U}(\gog)\otimes_{{\fam 2 U}(\gop)}\mV^\ch.
\end{eqnarray}

To summarize, the algebraic Fourier transform (the F-method) for generalized Verma modules allows to 
convert the explicit form of the action $(d\pi^\ch)(Z)$ to a system of partial differential 
equations $(d\tpi^\ch)(Z)$, i.e., it transforms the initial algebraic and combinatorial 
problem of computation of 
singular vectors in generalized Verma modules into analytic problem of solving a system of PDEs.
In particular, we introduce
\begin{align}
\label{eqn:sol}
 \mbox{\rm Sol}(\mathfrak{g},\mathfrak{g}';\mV^{\ch})
 :=\{f \in \Pol[\mathfrak {n}]\otimes \mV^\ch :
d\tilde\pi^\ch(Z) f = 0 \,\textrm{ for any } Z\in\gon'
\},
\end{align}
so the inverse Fourier transform gives
 an ${L'_{\mathbb R}}$-isomorphism
\begin{equation}
\label{eqn:phi}
\mbox{\rm Sol}(\mathfrak{g},\mathfrak{g}';\mV^{\ch})
 \overset \sim \longrightarrow
 M_{{\mathfrak {p}}}^{{\mathfrak {g}}}(\mV^{\ch})^{{\mathfrak {n}}'}.  
\end{equation}
An explicit form
of the action $d\tpi^\ch(Z)$ leads to 
a system  of differential equations for elements in $ \mbox{\rm Sol}$,
allowing to describe its structure completely
in some particular cases of interest (cf. \cite{kob, kope, koss}.)

In the dual language of differential operators acting on 
principal series representations, the set of $G_{\mathbb R}'$-equivariant differential operators 
$\Ind_{P_{\mathbb R}}^{G_{\mathbb R}}(\mV)\to\Ind_{P_{\mathbb R}'}^{G_{\mathbb R}'}(\mV')$ is in bijective
  correspondence with the space of all $(\gog', P_{\mathbb R}')$-homomorphisms
${\fam2 M}^{\gog'}_{\gop'}({\mV'}^\ch)\to{\fam2 M}^{\gog}_{\gop}(\mV^\ch).$

Now we describe the explicit choice of Lie algebras and representations
to which we apply the previous general procedure, and we also discuss the 
qualitative results on the underlying branching rule.
Let $n\in \mN$ be such that $n\geq 3$.
In the rest of the article $\gog$ denotes complexification of the real Lie algebra 
$so(n+1,1,\mR)$ of 
the connected simple real Lie group $G_{\mathbb R}=SO_o(n+1,1,\mR)$, and
$\gop$ its maximal parabolic subalgebra $\gop=\gol\ltimes\gon $, in the Dynkin 
diagrammatic notation for parabolic subalgebras given by omitting the first 
simple root of $\gog$. 
The Levi factor $\gol$ of $\gop$ is isomorphic to complexification of $so(n,\mR)\times\mR$ 
and the commutative 
nilradical $\gon$ (resp., the opposite nilradical $\gon_-$) is 
isomorphic to $\mC^n\simeq \mR^n\otimes_\mR\mC$. 

As for the matrix realization of this decomposition, 
we recall the Langlands decomposition 
$P_{\mathbb R}=L_{\mathbb R}N_{\mathbb R} = M_{\mathbb R}A_{\mathbb R}N_{\mathbb R}$.
The group $M_{\mathbb R}$ is isomorphic to $SO(n,{\mathbb R})$,
 and acts on ${\mathfrak {n}}_\mathbb{R} \simeq {\mathbb{R}}^n$
 by its fundamental vector representation
 preserving the quadratic form 
$\sum_{i=1}^{n}x_i^2$. The group $A_{\mR}$ is given by  
\begin{equation}
\label{eqn:grpA}
  A_{\mR}=\{
\left(
\begin{array}{ccc}
a& 0 & 0 \\
0 & I_n & 0\\
0 & 0 & a^{-1} 
\end{array}
\right) | \, a\in\mR^+\}\, \simeq \, \mR^+ .
\end{equation}
Then the elements $p\in P_{\mathbb R}$ are given by block triangular matrices
\begin{equation}
\label{eqn:pma}
  p=
\left(
\begin{array}{ccc}
a& \star & \star \\
0 & m & \star\\
0 & 0 & a^{-1} 
\end{array}
\right)
\end{equation}
with $a\in\mR^+,$ $m \in SO(n,{\mathbb R})$. 
Let $\{E_j\}_{j=1,\cdots,n}$ and $\{E_j^-\}_{j=1,\cdots,n}$
 be the standard basis of root vectors in ${\mathfrak {n}}_{\mathbb R}$
and ${{\mathfrak {n}}_{\mathbb R}}_-$, respectively, given by
\begin{equation}
\label{eqn:Epm}
 E_j =  
\left(
\begin{array}{ccc}
0 & e_j & 0 \\
0 & 0 & -e^t_j\\
0 & 0 & 0 
\end{array}
\right),\quad 
E_j^- =  
\left(
\begin{array}{ccc}
0 & 0 & 0 \\
e^t_j & 0 & 0\\
0 & -e_j & 0 
\end{array}
\right), \quad 1\leq j\leq n .
\end{equation}
We identify ${{\mathfrak {n}}_{\mathbb R}}_-$
 with ${\mathbb{R}}^n$ 
 and likewise ${\mathfrak {n}}_{\mathbb R}$
 with ${\mathbb{R}}^n$
 by using coordinates:
 $$
 \gon_{\mathbb R}\simeq \{Z:Z=(z_1,\ldots,z_n)\},
\quad
 {\gon_{\mathbb R}}_-\simeq\{X: X^t=(x_1,\ldots,x_n)\}.  
 $$ 
Then we have for elements in $N_{\mathbb R}$ and ${N_{\mathbb R}}_-$
\begin{eqnarray}\label{2matrices}
n=\exp Z=
\begin{pmatrix}
1&Z&-\frac{|Z|^2}{2}\\
0&\Id&- Z^t\\
0&0&1
\end{pmatrix}\in N_{\mathbb R} ,
\,\, \nonumber \\
x=\exp X=\begin{pmatrix}
1&0&0\\
X&\Id&0\\
-\frac{|X|^2}{2}&-X^t &1
\end{pmatrix}
\in {N_{\mathbb R}}_-,
\end{eqnarray}
where we set $|X|^2:=X^t X$  and $|Z|^2:=ZZ^t.$
By slight abuse of notation, we write $p = man$ for \eqref{eqn:pma}.

Then the case of our interest is given by pairs of complexified Lie algebras
\begin{eqnarray}
(\gog\oplus\gog,\gop\oplus\gop), \quad 
\mathrm{diag}(\gog,\gop)=(\mathrm{diag}(\gog),\mathrm{diag}(\gop)),
\end{eqnarray} 
where $\mathrm{diag}$ denotes the diagonal embedding, and the main task of 
the present article concerns the branching problem for the family 
of scalar generalized Verma 
${\fam2 U}(\gog\oplus\gog)$-modules induced from characters 
of $\gop\oplus\gop$. In particular, we shall 
classify (out of a discrete subset of inducing parameters) the scalar-valued 
singular vectors in scalar generalized 
Verma modules for orthogonal Lie algebras, i.e., the singular 
vectors transforming in the trivial representation of the simple part 
of Levi factor $\gol$.
 
Let $\mC_\lambda$ ($\mC_\mu$) denote the one-dimensional representation of $P_{\mathbb R}$ 
given by $p=man\mapsto a^\lambda$ ($p=man\mapsto a^\mu$), and 
$\chi_\lambda: \gop_{\mathbb R} \to {\mathbb C}_\lambda$ ($\chi_\mu: \gop_{\mathbb R} \to {\mathbb C}_\mu$)
its differential given by multiplication by $\lambda$ ($\mu$, respectively). An inducing character $\chi_{\lambda,\mu}$ 
of $\gop\oplus\gop$ is determined by two complex 
characters $\chi_\mu,\chi_\lambda$, $\lambda,\mu\in{\mathbb C}$,
\begin{eqnarray}
\chi_{\lambda,\mu}\equiv(\chi_\lambda,\chi_\mu) : \gop\oplus\gop & \rightarrow & \mC ,
\\ \nonumber
(p_1,p_2)& \mapsto & \chi_\lambda(p_1)\otimes\chi_\mu(p_2)\in \mbox{End}(\mC_\lambda\otimes\mC_\mu),
\end{eqnarray}
and the generalized Verma ${\fam2 U}(\gog\oplus\gog)$-module induced from
character $(\chi_\lambda,\chi_\mu)$ is
\begin{eqnarray}
{\fam2 M}^{\gog\oplus\gog}_{\gop\oplus\gop}(\mC_\lambda\otimes\mC_\mu)
\equiv{\fam2 M}_{\lambda,\mu}(\gog\oplus\gog,\gop\oplus\gop)=
{\fam2 U}(\gog\oplus\gog)\otimes_{{\fam2 U}(\gop\oplus\gop)}(\mC_\lambda\otimes\mC_\mu),
\end{eqnarray}
where $\mC_\lambda\otimes\mC_\mu$ is a $1$-dimensional representation $(\chi_\lambda,\chi_\mu)$
of $\gop\oplus\gop$. 
As a vector space, ${\fam2 M}_{\lambda,\mu}(\gog\oplus\gog,\gop\oplus\gop)$ is isomorphic to 
the symmetric algebra $S(\gon_-\oplus\gon_-)$, where 
$\gon_-\oplus\gon_-$ is the complement of $\gop\oplus\gop$ in $\gog\oplus\gog$.
Notice that we have an isomorphism
$$
(\gon_-\oplus\gon_-)/((\gon_-\oplus\gon_-)\cap \mathrm{diag}(\gog))\simeq \gon_-
$$
of $\mathrm{diag}(\gol)$-quotient modules. 

The symmetric algebra $S\big((\gon_-\oplus\gon_-)/((\gon_-\oplus\gon_-)\cap \mathrm{diag}(\gog))\big)$ 
decomposes as $\mathrm{diag}(\gol)$-module on irreducible submodules, with higher multiplicities in general. 
In particular, each 
$\mathrm{diag}(\gol)$-module realized in homogeneity $k$ polynomials also appears in polynomials of homogeneity 
$(k+2)$, $k\in\mN$.
As we have already explained, we focus on the case of $1$-dimensional representations
$\mV_{\lambda,\mu}\simeq \mC_\lambda\otimes\mC_\mu$ regarded as $(\gop\oplus\gop)$-modules
(with the trivial action of the simple part of Levi subalgebra and the nilradical) and 
$\mV_{\lambda'}\simeq\mC_\nu$ as $\mathrm{diag}(\gol)$-modules ($\lambda,\mu,\nu\in\mC$),
and it is a result in classical invariant theory (see \cite{gw}, \cite{kp}) that 
for each even homogeneity there is just
one $1$-dimensional module. Because $\gon_-$ is as $(\mathrm{diag}(\gol)/[\mathrm{diag}(\gol),\mathrm{diag}(\gol)])$-module isomorphic 
to the character $\mC_{-1}$, 
the following holds true in the Grothendieck group of ${\fam2 O}^{\gop}$, $\gop\simeq \mathrm{diag}(\gop)$. 
As a consequence of \cite{ko}, Theorem $3.9$, we have
\begin{corollary}
\label{charform}
Let 
$$\gog\oplus\gog = so(n+2,\mC)\oplus so(n+2,\mC),\quad {\mathrm{diag}}(\gog)\simeq so(n+2,\mC)
$$ 
with standard maximal parabolic subalgebras $\gop\oplus\gop,\, \mathrm{diag}(\gop)$ 
given by omitting the first simple root in the corresponding Dynkin diagram. 

Then
the multiplicity $m(\nu,(\lambda,\mu))$ of ${\fam2 M}^{\gog}_{\gop}(\mC_\nu)\equiv{\fam2 M}_{\nu} (\gog,\gop)$ 
in ${\fam2 M}_{\lambda,\mu} (\gog\oplus\gog,\gop\oplus\gop)$ is equal to 
one for $\nu=\lambda+\mu -2j,\, j\in{\mathbb N}$, and zero otherwise.
In the Grothendieck group $K({\fam2 O}^\gop)$ of the Bernstein-Gelfand-Gelfand parabolic category 
${\fam2 O}^{\gop}$ holds
\begin{eqnarray}\label{diag-branch-tk}
{\fam2 M}_{\lambda,\mu} (\gog\oplus\gog,\gop\oplus\gop)|_{\mathrm{diag}(\gog)}
\simeq \bigoplus_{j\in\mN}
{\fam2 M}_{\lambda+\mu -2j} (\gog,\gop).
\end{eqnarray}
\end{corollary}

 Although we work in one specific signature $(n+1,1)$, the results are easily extended 
to any real form of arbitrary signature. 


\section{The construction of singular vectors for diagonal branching rules applied to scalar 
generalized Verma modules for $so(n+1,1,\mR)$} 

The rest of the article is devoted to the construction of a class of 
scalar valued singular vectors, whose abstract existence was concluded in Section $2$, 
Corollary \ref{charform}, using the tool of algebraic Fourier transform
(F-method) for generalized Verma modules reviewed in Section $2$.
This can be regarded as a quantitative part of our
diagonal branching problem.

\subsection{Description of the representation}

In this subsection we describe the representation of 
$\gog\oplus\gog$ on scalar generalized Verma modules 
\begin{eqnarray}\label{tensorproduct}
{\fam2 M}_{\lambda,\mu}(\gog\oplus\gog,\gop\oplus\gop)=
{\fam2 U}(\gog\oplus\gog)\otimes_{{\fam2 U}(\gop\oplus\gop)}\mC_{\lambda,\mu}
\simeq{\fam2 M}_\lambda(\gog,\gop)\otimes{\fam2 M}_\mu(\gog,\gop)
\end{eqnarray}
for $\mC_{\lambda,\mu}=\mC_\lambda\otimes\mC_\mu$ in its Fourier image, i.e., we apply the 
F-method explained in Section $2$.
The first goal is to describe the action by elements in the nilradical $\mathrm{diag}(\gon)$ of
$\mathrm{diag}(\gop)$ in terms of
differential operators acting on the Fourier image of ${\fam2 M}_{\lambda,\mu}(\gog\oplus\gog,\gop\oplus\gop)$.
It can be derived from the explicit form of the action on the induced 
representation realized in the non-compact picture, and it follows from ($\ref{tensorproduct}$)
that the problem can be reduced to the question on each component of the tensor 
product separately.
Let us consider the complex representation $\pi_\la$ of $G_{\mathbb R}=SO_o(n+1,1,\mR)$ 
on $\Ind_{P_{\mathbb R}}^{G_{\mathbb R}}(\mC_\lambda)$, $\la\in \mC$, induced from
the character $p=man\mapsto a^\lambda$, $p\in P_{\mathbb R}$, on one dimensional representation 
space $\mC_\lambda\simeq\mC$. Here $a\in A_{\mathbb R}$
is the abelian subgroup in the Langlands decomposition $P_{\mathbb R}=M_{\mathbb R}A_{\mathbb R}N_{\mathbb R}$, 
$M_{\mathbb R}\simeq SO(n,\mR),\, N_{\mathbb R}\simeq\mR^n$. The character of $P_{\mathbb R}$
is trivial on $M_{\mathbb R}$ and $N_{\mathbb R}$, and its value on $a\in A_{\mathbb R}\simeq {\mathbb R}^+$ is $a^\lambda$
(i.e., it is a complex character of $A_{\mathbb R}$ on ${\mathbb C}_\lambda$.)  

Let $x_j$, $j=1,\dots ,n$, be the coordinates with respect to the standard basis $\{E_j^-\}_{j=1,...,n}$ of root 
vectors on $\gon_-$, and $\xi_j$, $j=1,\dots ,n$, the coordinates on the Fourier transform of $\gon_-$. 
We consider the family of differential operators 
\begin{eqnarray}
Q_j(\la)=-\frac{1}{2} |x|^2\partial_j+x_j(-\lambda +\sum_kx_k\partial_k),\, j=1,\ldots, n, 
\end{eqnarray}
\begin{equation}
\label{dpti}
P_j^\xi(\lambda) = {\mathbf i}\left(\frac{1}{2}\xi_j\triangle^\xi + (\lambda-\mE^\xi)\partial_{\xi_j}\right),\, j=1,\dots ,n,
\end{equation}
where $|x|^2=x_1^2+\dots +x_n^2$,
 $$\triangle^\xi =\partial^2_{\xi_1}+\dots +\partial^2_{\xi_{n}}$$
is the Laplace operator of positive signature, $\partial_j=\frac{\partial}{\partial x_j}$ 
and  $\mE^\xi=\sum_k\xi_k\partial_{\xi_k}$ is 
the Euler homogeneity operator (${\mathbf i}\in\mC$ denotes the complex unit.) 
Via the exponential map, the non-compact picture of the induced representation 
$\Ind_{P_{\mathbb R}}^{G_{\mathbb R}}(\mC_\lambda)$ is given by
$$
C^\infty({\gon_{\mathbb R}}_-,\mC_{\lambda}) \simeq C^\infty({\gon_{\mathbb R}}_-,\mC)\otimes \mC_{\lambda}.
$$
The following result is a routine computation, cf. \cite{koss}: 
 
\begin{lemma}
Let us denote by $E_j\in\gon_{\mathbb R}$ the standard basis elements, 
$j=1,\dots ,n$. Then 
$E_j$ act on $C^\infty({\gon_{\mathbb R}}_-,\mC)\otimes \mC_{\lambda}$ by 
\begin{eqnarray}\label{dpio}
d\pi_{\lambda}(E_j)(s\otimes v)=Q_j(\la)(s)\otimes v, \,
s\in{C^\infty}({\gon_{\mathbb R}}_-,\mC),\, {v}\in \mC_{\lambda},  
\end{eqnarray}
and the action of $d\tilde{\pi}^\vee_\lambda$
on $\Pol[\xi_1,\ldots,\xi_n]\otimes \mC^\vee_{\la}$ is given by
\begin{eqnarray}\label{dtpio}
d\tilde{\pi}^\vee_{\la}(E_j)(f\otimes u)=P_j^\xi(\la)(f)\otimes u,\,
f\in\Pol[\xi_1,\ldots,\xi_n],\,
u\in \mC_{\lambda}^\ch.  
\end{eqnarray}
\end{lemma}
As for the action of remaining basis elements of $\gog$ in the Fourier image of 
the induced representation, 
the action of $\gon_-$ is given by multiplication by coordinate functions, the standard basis 
elements of the simple part of the Levi factor $\gol^s=[\gol,\gol]\simeq so(n,\mC)$ realized by
matrices 
$$
 M_{ij} =  
\left(
\begin{array}{ccc}
0 & 0 & 0 \\
0 & {(M_{ij})}_{rs}=\delta_{ir}\delta_{js}-\delta_{is}\delta_{jr} & 0   \\
0 & 0 & 0 
\end{array}
\right),\quad i,j,r,s=1,...,n, i<j,
$$
act by first order differential operators 
$$
M^\xi_{ij}=(\xi_j\partial_{\xi_i}-\xi_i\partial_{\xi_j}),\, i,j=1,\dots ,n
$$ 
and the basis element of the Lie algebra of $A_{\mathbb R}$ given by the matrix
$$
 E =  
\left(
\begin{array}{ccc}
1 & 0 & 0 \\
0 & 0 & 0   \\
0 & 0 & -1 
\end{array}
\right)
$$
acts as the homogeneity
operator, $\mE^\xi=\sum_{i=1}^n\xi_i\partial_{\xi_i}$.   

We introduce the diagonal embedding of the Lie algebra $\gog_{\mathbb R}$ as the differential of 
the diagonal map for the Lie group $G_{\mathbb{R}}$
$$
G_{\mathbb{R}}\hookrightarrow G_{\mathbb{R}}\times G_{\mathbb{R}},\quad g\mapsto (g,g),
$$
given by
$$
\mathrm{diag}(\gog_\mR): X\mapsto X\otimes 1 + 1\otimes X,\quad \mbox{for all}\quad X\in\gog_\mR 
$$
and termed diagonaly embedded Lie algebra of $\gog_\mR$. 
This notion restricts to any Lie subalgebra of $\gog_\mR$.
Then the action of $\mathrm{diag}(\gog_\mR)$ in the representation 
$d\tilde{\pi}^\vee_{\la}\otimes d\tilde{\pi}^\vee_{\mu}$ on 
$\Pol[\xi_1,\ldots,\xi_n]\otimes \mC^\vee_{\la}\otimes \Pol[\nu_1,\ldots,\nu_n]\otimes \mC^\vee_{\mu}$,
the algebra of complex valued polynomials on ${\gon}\oplus{\gon}$ with coordinates 
$\xi_i$ resp. $\nu_i$ on 
the first resp. second copy of ${\gon}$ in ${\gon}\oplus{\gon}$, by
particular generators of the diagonal subalgebra $\mathrm{diag}(\gog_{\mathbb R})$ is given by 
\begin{enumerate}
\item
The multiplication by 
\begin{eqnarray}
(\xi_j\otimes 1)+(1\otimes\nu_j),\, j=1,\dots ,n\, ,
\end{eqnarray}  
for the elements of $\mathrm{diag}({\gon_\mR}_-)$,
\item
The action by first order differential operators with linear coefficients
\begin{eqnarray}
 M^{\xi ,\nu}_{ij} &=& (M^\xi_{ij}\otimes 1)+(1\otimes M^\nu_{ij})
\nonumber \\
& = &(\xi_j\partial_{\xi_i}-\xi_i\partial_{\xi_j})\otimes 1+
1\otimes (\nu_j\partial_{\nu_i}-\nu_i\partial_{\nu_j})\, ,
\end{eqnarray}  
$i,j=1,\dots ,n$, for the elements of the simple part of $\mathrm{diag}(\gol_\mR)$ and
\begin{eqnarray}
\mE^\xi\otimes 1+1\otimes\mE^\nu=\sum_{i=1}^n(\xi_i\partial_{\xi_i}\otimes 1+1\otimes \nu_i\partial_{\nu_i})
\end{eqnarray}
for the generator of $\mathrm{diag}(\gol_\mR)/[\mathrm{diag}(\gol_\mR),\mathrm{diag}(\gol_\mR)]$,
\item
The action by second order differential operators with at most linear coefficients
\begin{eqnarray}
P_j^{\xi, \nu}(\lambda ,\mu)&=&(P_j^\xi(\lambda)\otimes 1)+(1\otimes P_j^\nu(\mu))
\nonumber \\
&=&
{\mathbf i}(\frac{1}{2}\xi_j\triangle^\xi+(\lambda-\mE^\xi)\partial_{\xi_j})\otimes 1
\nonumber \\
& & +{\mathbf i}1\otimes(\frac{1}{2}\nu_j\triangle^\nu+(\mu-\mE^\nu)\partial_{\nu_j}),    
\end{eqnarray}
$j=1,\dots ,n$ for the elements $\mathrm{diag}(\gon_\mR)$. 
\end{enumerate}


\subsection{Reduction to scalar differential equation in two variables}

It follows from the discussion in Section $2$ that $\mathrm{diag}(\gol)$-modules 
for the diagonal branching rules, inducing singular vectors for generalized 
scalar Verma modules, are one dimensional. This means that
they are annihilated by $\mathrm{diag}(\gol^s)\simeq so(n,\mC)$, the simple part of the diagonal 
Levi subalgebra
$\mathrm{diag}(\gol)\simeq so(n,\mC)\times\mC$. It follows 
that the singular vectors are invariants of $\mathrm{diag}(\gol^s)$ acting 
diagonally on the algebra of polynomials on $\gon\oplus\gon$ regarded 
as a $\gol\oplus\gol$-module. 
The following result is a special case of the first fundamental theorem in classical invariant 
theory, see e.g. \cite{gw}, \cite{kp}.   
\begin{lemma}
Let $(V,\langle,\rangle)$ be a finite dimensional complex vector space with bilinear form 
$\langle,\rangle$ and $SO(V)$ the Lie 
group of automorphisms of $(V,\langle,\rangle)$. Then the subalgebra $\Pol[V\oplus V]^{SO(V)}$ of
$SO(V)$-invariants in the complex polynomial algebra 
$\Pol[V\oplus V]$ (with $SO(V)$ acting diagonally on $V\oplus V$) is
polynomial algebra generated by $\langle\xi,\xi\rangle$, $\langle\xi,\nu\rangle$ and
$\langle\nu,\nu\rangle$. Here we use the convention that $\xi$ is a vector in the
first summand $V$ of the direct sum $V\oplus V$ and $\nu$ in the second summand.
Therefore, there is an algebra isomorphism
$$
\Pol[V\oplus V]^{SO(V)}\simeq \Pol[\langle\xi,\nu\rangle, \langle\xi,\xi\rangle, \langle\nu,\nu\rangle].
$$
\end{lemma}
In our case, the complex polynomial algebra is $\Pol[\xi_1,\dots ,\xi_n,\nu_1,\dots ,\nu_n]$
and we use the notation $\Pol[r,s,t]$ for the subalgebra of invariants generated by
\begin{eqnarray}
r:=\langle\xi,\nu\rangle=\sum_{i=1}^n\xi_i\nu_i,\, \nonumber \\
s:=\langle\xi,\xi\rangle=\sum_{i=1}^n\xi_i\xi_i,\, \nonumber \\
t:=\langle\nu,\nu\rangle=\sum_{i=1}^n\nu_i\nu_i.
\end{eqnarray}
The task of the present subsection is to rewrite the operators 
$P^{\xi,\nu}_j(\lambda,\mu)$ in the variables $r,s,t$, i.e., to
reduce the action of $P^{\xi,\nu}_j(\lambda,\mu)$ from the polynomial 
ring to the algebra of $\mathrm{diag}(\gol^s)$-invariants on $\gon\oplus\gon$.  

We compute 
\begin{eqnarray}
 \partial_{\nu_i}r=\xi_i,\, \partial_{\xi_i}r=\nu_i,
 \partial_{\nu_i}s=0,\, \partial_{\xi_i}s=2\xi_i,
 \partial_{\nu_i}t=2\nu_i,\, \partial_{\xi_i}t=0,
\end{eqnarray}
and 
\begin{eqnarray}
\partial_{\xi_i}=\nu_i{\partial_r}+2\xi_i{\partial_s},
\triangle^\xi=t\partial_r^2+4r\partial_r\partial_s+2n\partial_s+4s\partial_s^2,\, i=1,\dots ,n.
\end{eqnarray}
Note that analogous formulas for $\partial_{\nu_i}$, $\triangle^\nu$ can be 
obtained from those for $\xi$ by applying the change of variables 
\begin{eqnarray} 
\xi_i\longleftrightarrow \nu_i,\, s\longleftrightarrow t,\, r\longleftrightarrow r.
\end{eqnarray}
We also have for all $i=1,\dots ,n$
\begin{eqnarray}
\mE^\xi\partial_{\xi_i}=\nu_i(\mE^r+2\mE^s)\partial_r+\xi_i(2\mE^r+4\mE^s+2)\partial_s,
\end{eqnarray}
and so taking all together we arrive at the operators
\begin{eqnarray}\label{initialFsystem}
& & P^{r,s,t}_i(\lambda,\mu)=\xi_i(\frac{1}{2}t\partial_r^2+
(n+2\lambda-2-2\mE^s)\partial_s-(\mE^r+2\mE^t-\mu)\partial_r)
\nonumber \\
& & +\nu_i(\frac{1}{2}s\partial_r^2+(n+2\mu-2-2\mE^t)\partial_t
-(\mE^r+2\mE^s-\lambda)\partial_r)
\end{eqnarray}
acting on complex polynomial algebra $\Pol[r,s,t]$, $i=1,\dots ,n$. A suitable
linear combination of this vector-valued system of equations ($i=1,\dots ,n$) 
leads to operators
\begin{eqnarray}\label{formersystem}
& & P^{r,s,t}_{\xi}(\lambda,\mu):=\sum_{i=1}^n\xi_iP^{r,s,t}_i(\lambda,\mu)=
s(\frac{1}{2}t\partial_r^2+(n+2\lambda-2-2\mE^s)\partial_s
\nonumber \\
& & -(\mE^r+2\mE^t-\mu)\partial_r)
+r(\frac{1}{2}s\partial_r^2+(n+2\mu-2-2\mE^t)\partial_t
-(\mE^r+2\mE^s-\lambda)\partial_r),
\nonumber \\
& & P^{r,s,t}_{\nu}(\lambda,\mu):=\sum_{i=1}^n\nu_iP^{r,s,t}_i(\lambda,\mu)=
r(\frac{1}{2}t\partial_r^2+(n+2\lambda-2-2\mE^s)\partial_s
\nonumber \\
& &-(\mE^r+2\mE^t-\mu)\partial_r)+t(\frac{1}{2}s\partial_r^2+(n+2\mu-2-2\mE^t)\partial_t
-(\mE^r+2\mE^s-\lambda)\partial_r).
\nonumber \\ \label{system}
\end{eqnarray} 
Notice that the second equation follows from the first one by the action 
of involution 
$$
\lambda\longleftrightarrow\mu,\, s\longleftrightarrow t,\, r\longleftrightarrow r.
$$
In what follows we construct a set of homogeneous polynomial solutions of 
$P^{r,s,t}_{\xi}(\lambda,\mu)$, $P^{r,s,t}_{\nu}(\lambda,\mu)$ solving the
system $\{P^{r,s,t}_i(\lambda,\mu)\}_i$, $i=1,\dots ,n$. The uniqueness of 
the solution for the generic values of the inducing parameters 
implies the unique solution of the initial system of PDEs \eqref{initialFsystem}.
Notice that (\ref{system}) is the system of differential operators 
preserving the 
space of homogeneous polynomials in the variables $r,s,t$, i.e. 
$P^{r,s,t}_{\xi}(\lambda,\mu)$, $P^{r,s,t}_{\nu}(\lambda,\mu)$ commute 
with $\mE^{r,s,t}:=\mE^r+\mE^s+\mE^t$. 


\subsection{Polynomial solutions of the differential equation in two variables produced by the F-method}

We start with a couple of simple examples.
\begin{example}\label{linear-example}
Let us consider the polynomial of homogeneity one, 
$$p(r,s,t)=Ar+Bs+Ct, \,\,A,B,C\in\mC .$$
The application of $P_i^{r,s,t}(\lambda,\mu)$  yields
\begin{eqnarray}
P_i^{r,s,t}(\lambda,\mu)(Ar+Bs+Ct)=\xi_i(B(n+2\lambda -2)+A\mu)
+\nu_i(C(n+2\mu-2)+A\lambda)
\end{eqnarray}
for all $i=1,\dots ,n$. When $A$ is normalized to be equal to $1$, we get
$$
C=-\frac{\lambda}{n+2\mu-2}, \, B=-\frac{\mu}{n+2\lambda-2}.
$$ 
If $n+2\mu-2\not= 0$ and $n+2\lambda-2\not= 0$, there is a unique homogeneous 
solution of $P^{r,s,t}_i(\lambda,\mu)p(r,s,t)=0$ for all $i=1,\dots ,n$  given 
by
\begin{eqnarray}
p(r,s,t)=(n+2\lambda-2)(n+2\mu-2)r
-{\mu}(n+2\mu-2)s-{\lambda}(n+2\lambda-2)t.
\end{eqnarray}
\end{example}
\begin{example}\label{quadratic-example}
Let 
$$p(r,s,t)=Ar^2+Bs^2+Ct^2+Drs+Est+Frt, \,\, A,B,C,D,E,F\in\mC$$ 
be a general polynomial of homogeneity two. We have 
\begin{eqnarray}
& & P^{r,s,t}_i(\lambda,\mu)p(r,s,t)=\xi_i[r(D(n+2\lambda-2)+A2(\mu-1))
+s(B2(n+2\lambda-4)+D\mu)
\nonumber \\
& & +t(A+E(n+2\lambda-2)+F(\mu-2))],
\nonumber \\
& &+\nu_i[r(F(n+2\mu-2)+A2(\lambda-1))
+s(A+E(n+2\mu-2)+D(\lambda-2))
\nonumber \\
& & +t(C2(n+2\mu-4)+F\lambda)],
\end{eqnarray}
for all $i=1,\dots ,n$. The equations 
$$\sum_i\xi_iP^{r,s,t}_{\xi_i}(\lambda,\mu)=0,\,\,
\sum_i\nu_iP^{r,s,t}_{\nu_i}(\lambda,\mu)=0
$$
are equivalent to two systems of linear equations:
\begin{eqnarray}
& & D(n+2\lambda-2)+A2(\mu-1)+A+E(n+2\mu-2)+D(\lambda-2)=0,
\nonumber \\
& & F(n+2\mu-2)+A2(\lambda-1)=0,
\nonumber \\
& & B2(n+2\mu-4)+D\mu=0,
\nonumber \\
& & A+E(n+2\lambda-2)+F(\mu-2)=0,
\nonumber \\
& & C2(n+2\mu-4)+F\lambda=0,
\end{eqnarray}
resp.
\begin{eqnarray}
& & D(n+2\lambda-2)+A2(\mu-1)=0,
\nonumber \\
& & B2(n+2\mu-4)+D\mu=0,
\nonumber \\
& & A+E(n+2\lambda-2)+F(\mu-2)+F(n+2\mu-2)+A2(\lambda-1)=0,
\nonumber \\
& & A+E(n+2\mu-2)+D(\lambda-2)=0,
\nonumber \\
& & C2(n+2\mu-4)+F\lambda=0.
\end{eqnarray}
Both systems are equivalent under the involution 
$$
A\longleftrightarrow A,\, E\longleftrightarrow E,\, D\longleftrightarrow F,\,
B\longleftrightarrow C,\, \lambda\longleftrightarrow \mu ,
$$
and for $n+2\mu-2\not= 0, n+2\lambda-2\not= 0, n+2\mu-4\not= 0$ and $n+2\lambda-4\not= 0$
there is a unique solution invariant under this involution
\begin{eqnarray} 
& & A=1,\, F=\frac{-2(\lambda -1)}{n+2\mu-2},\, C=\frac{\lambda(\lambda -1)}{(n+2\mu-2)(n+2\mu-4)},
\nonumber \\ 
& & E=2\frac{(\lambda -2)(\mu-2)-(1+\frac{n}{2})}{(n+2\mu-2)(n+2\lambda-2)},\,
D=\frac{-2(\mu -1)}{(n+2\lambda-2)},\, 
\nonumber \\
& & B=\frac{\mu(\mu -1)}{(n+2\lambda-2)(n+2\lambda-4)}.
\end{eqnarray}
To conclude, in the case when $n+2\mu-2\not= 0, n+2\lambda-2\not= 0, n+2\mu-4\not= 0$ 
and $n+2\lambda-4\not= 0$, the polynomial 
\begin{eqnarray}
p(r,s,t)&=&(n+2\lambda-2)(n+2\lambda-4)(n+2\mu-2)(n+2\mu-4)r^2
\nonumber \\
& & +\mu(\mu-1)(n+2\mu-2)(n+2\mu-4)s^2
\nonumber \\
& & +\lambda(\lambda-1)(n+2\lambda-2)(n+2\lambda-4)t^2
\nonumber \\
& & -2(\mu-1)(n+2\lambda-4)(n+2\mu-2)(n+2\mu-4)rs
\nonumber \\
& & +2((\lambda -2)(\mu-2)-(1+\frac{n}{2}))(n+2\lambda-4)(n+2\mu-4)st
\nonumber \\
& &  -2(\lambda-1)(n+2\lambda-2)(n+2\lambda-4)(n+2\mu-4)rt
\end{eqnarray}
is the unique solution of $P^{r,s,t}_i(\lambda,\mu)p(r,s,t)=0$ 
of homogeneity two.
\end{example}

We now return back to the case of general homogeneity. Let 
\begin{eqnarray}\label{non-homsol}
p=p(r,s,t)=\sum_{0\leq i,j,k\leq N|i+j+k=N}A_{i,j}s^it^jr^k
\end{eqnarray}
be a homogeneous polynomial of degree $N$, $deg(p)=N$, and write 
\begin{eqnarray}\label{subuv}
& & p=r^Np(1,\frac{s}{r},\frac{t}{r})=r^N\tilde{p}(u,v),\,\, u:=\frac{s}{r}, v:=\frac{t}{r},
\nonumber \\ 
& & \tilde{p}(u,v)=\sum_{0\leq i,j|i+j\leq N}A_{i,j}u^iv^j,
\end{eqnarray}
where $\tilde{p}$ is a polynomial of degree $N$. 
The dehomogenesation $(r,s,t)\to (r,u,v)$ is governed by coordinate change
\begin{eqnarray}
u:=\frac{s}{r},\, v:=\frac{t}{r},\, r:=r,
\end{eqnarray}
such that
\begin{eqnarray}
& & \partial_s\to\frac{1}{r}\partial_u,\, \partial_t\to\frac{1}{r}\partial_v,\,
\partial_r\to-\frac{1}{r}u\partial_u-\frac{1}{r}v\partial_v+\partial_r\, ,  
\nonumber \\
& & \mE^s\to \mE^u,\, \mE^t\to \mE^v,\, \mE^r\to -\mE^u-\mE^v+\mE^r\, .
\end{eqnarray}
Then the particular terms in $P^{r,s,t}_\xi(\lambda ,\mu)$, see \eqref{formersystem}, transform as 
\begin{eqnarray}
\frac{1}{2}st\partial_r^2 &\to &\frac{1}{2}uv(\mE^u+\mE^v-\mE^r+1)(\mE^u+\mE^v-\mE^r),
\nonumber \\
 -r(\mE^r+2\mE^s-\lambda)\partial_r &\to& (\mE^u+\mE^v-\mE^r)(\mE^u-\mE^v+\mE^r-\lambda-1),
\nonumber \\
 r(n+2\mu-2-2\mE^t)\partial_t &\to& (n+2\mu-2-2\mE^v)\partial_v,
\nonumber \\
 \frac{1}{2}rs\partial_r^2 &\to& \frac{1}{2}u(\mE^u+\mE^v-\mE^r+1)(\mE^u+\mE^v-\mE^r),
\nonumber \\
 -s(\mE^r+2\mE^t-\mu)\partial_r &\to& u(\mE^u+\mE^v-\mE^r)(-\mE^u+\mE^v+\mE^r-\mu-1),
\nonumber \\
 s(n+2\lambda-2-2\mE^s)\partial_s &\to& (n+2\lambda-2\mE^u)\mE_u,
\end{eqnarray}
and when acting on a polynomial of homogeneity $N$, $p(r,s,t)=r^N{\tilde p}(u,v)$ for a 
polynomial ${\tilde p}(u,v)$ of degree $N$ in $u,v$, $\mE^r=N$ and we get
\begin{eqnarray}\label{conversion}
 P^{u,v}_\xi (\lambda,\mu) &=& \frac{1}{2}uv(\mE^u+\mE^v-N+1)(\mE^u+\mE^v-N)
\nonumber \\
& & -({\mE^u})^2+\mE^u(n+\lambda-1)+(\mE^v-N)(-\mE^v+N-\lambda-1)
\nonumber \\
& & +(n+2\mu-2-2\mE^v)\partial_v
\nonumber \\
& & +\frac{1}{2}u(\mE^u+\mE^v-N)(-\mE^u+3\mE^v+N-2\mu-1).
\end{eqnarray}
Similarly, we have
\begin{eqnarray}
 P^{u,v}_\nu (\lambda,\mu) &=& \frac{1}{2}uv(\mE^u+\mE^v-N+1)(\mE^u+\mE^v-N)
\nonumber \\
& & -({\mE^v})^2+\mE^v(n+\mu-1)+(\mE^u-N)(-\mE^u+N-\mu-1)
\nonumber \\
& & +(n+2\lambda-2-2\mE^u)\partial_u
\nonumber \\
& & +\frac{1}{2}v(\mE^v+\mE^u-N)(-\mE^v+3\mE^u+N-2\lambda-1).
\end{eqnarray}
Let us denote $A_{i,j}(\lambda ,\mu)$ the coefficient by monomial $u^iv^j$ 
in the polynomial $\tilde{p}(u,v)$. 

The assumption $A_{i,j}(\lambda ,\mu)=A_{j,i}(\mu,\lambda)$, 
combined with the symmetry between $P^{u,v}_\xi (\lambda,\mu)$ and 
$P^{u,v}_\nu (\lambda,\mu)$, allows to restrict to the action of 
$P^{u,v}_\xi (\lambda,\mu)$ on a polynomial of degree $N$ of the form
\begin{eqnarray}\label{formsingvect}
{\tilde p}(u,v)=\sum_{i,j|0\leq |i+j|\leq N}A_{i,j}(\lambda,\mu)u^iv^j,\,\, 
A_{i,j}(\lambda ,\mu)=A_{j,i}(\mu,\lambda) .
\end{eqnarray}
Consequently, we convert the differential equation (\ref{conversion}) into the four-term 
functional relation 
\begin{eqnarray}\label{rfe}
& & [\frac{1}{2}(i+j-N-1)(i+j-N-2)]A_{i-1,j-1}(\lambda ,\mu)
\nonumber \\
& & +[(-i^2+i(n+\lambda-1)+(j-N)(-j+N-\lambda-1)]A_{i,j}(\lambda ,\mu)
\nonumber \\
& & +[(j+1)(n+2\mu-2-2j)]A_{i,j+1}(\lambda ,\mu)
\nonumber \\
& & +[\frac{1}{2}(i+j-N-1)(-i+3j+N-2\mu)]A_{i-1,j}(\lambda ,\mu)=0
\end{eqnarray}
for $i,j=1,\dots ,N$ and $j\geq i$, which recursively computes $A_{i,j+1}(\lambda ,\mu)$
in terms of $A_{i-1,j-1}(\lambda ,\mu),A_{i-1,j}(\lambda ,\mu)$ and $A_{i,j}(\lambda ,\mu)$.

As for the normalization of $A_{i,j}(\lambda ,\mu)$, a singular 
vector can be normalized by multiplication by common denominator 
resulting in the coefficients valued in $\mathrm{Pol}[\lambda,\mu]$ rather than its 
quotient field $\mC(\lambda,\mu)$. 

As we shall prove in the next Theorem, a  
consequence of (\ref{rfe}) is the uniqueness of 
its solution in the range $\lambda,\mu\in\mC\setminus\{m-\frac{n}{2}\, |\, m\in\mN\}$. 
We observe that the uniqueness of solution fails for $\lambda,\mu\in\{m-\frac{n}{2}\, |\, m\in\mN\}$,
which indicates the appearance of a non-trivial composition structure in the 
branching problem for generalized Verma modules.     
In the following Theorem we construct a set of singular vectors, which 
will be the representatives realizing abstract character 
formulas of the diagonal branching problem in Corollary \ref{charform}.
\begin{theorem}\label{A-comp}
Let us assume $\lambda,\mu\in\mC\setminus\{m-\frac{n}{2}\, |\, m\in\mN\}$, $N\in\mN$, 
and introduce the Pochhammer symbol $(x)_l=x(x+1)\dots (x+l-1),\,l\in\mN$, for $x\in\mC$. 
The four-term functional equation (\ref{rfe}) for the
coefficients $\{A_{i,j}(\lambda,\mu)\}_{i,j\in\{1,\dots ,N\}}$ 
of (scalar valued) singular vectors fulfilling 
$$
A_{j,i}(\lambda,\mu)=A_{i,j}(\mu,\lambda),\, j\geq i,
$$   
has the unique non-trivial solution given by the formula  
\begin{eqnarray} \label{solut} 
& & A_{i,j}(\lambda ,\mu)=  
\nonumber \\
& & \frac{\Gamma(i+j-N)\Gamma(1-\frac{n}{2}-\mu)\Gamma(1-i+j-N+\lambda)\Gamma(\lambda+\frac{n}{2}-i)}
{2^{i+j}(-1)^{i+j}i!j!\,\Gamma(-N)\Gamma(1-N+\lambda)\Gamma(1+j-\frac{n}{2}-\mu)\Gamma(\lambda+\frac{n}{2})}\cdot
\nonumber \\  
& & \sum_{k=0}^i(-1)^k {i\choose k}(j-i+1+k)_{i-k}(\lambda+\frac{n}{2}-i)_{i-k}(\mu-N+1)_k(\lambda-N+1-k)_k.  
\nonumber \\
\end{eqnarray}  
 \end{theorem}  
{\bf Proof:} 

Let us first discuss the uniqueness of the solution. The 
knowledge of $A_{i,j}(\lambda,\mu)$ for $i+j\leq k_0$ allows to compute
the coefficient $A_{i,j+1}(\lambda,\mu)$ with $i+j=k_0+1$ from the 
recursive functional equation, because of assumption 
$\lambda,\mu\notin\{m-\frac{n}{2}\, |\, m\in\mN\}$. The symmetry condition
for $A_{i,j}(\lambda,\mu)$ gives $A_{j,i}(\mu,\lambda)=A_{i,j}(\lambda,\mu)$
and the induction proceeds by passing to the computation of 
$A_{i,j+2}(\lambda,\mu)$. Note that all coefficients are proportional to 
$A_{0,0}(\lambda,\mu)$ and its choice affects their explicit form. 

The proof of the explicit form for $A_{i,j}(\lambda,\mu)$ is based on 
the verification of the recursion functional 
equation (\ref{rfe}). To prove that the left hand side of (\ref{rfe}) is trivial
is equivalent to the following check: up to a product of linear factors
coming from $\Gamma$-functions, the left hand side is the sum of 
four polynomials in $\lambda, \mu$. A simple criterion for the triviality
of a polynomial of degree $d$ we use is that it has $d$ roots (counted with multiplicity)
and the leading monomial in a corresponding variable has coefficient zero.

It is straightforward but tedious to check that the left hand side 
of (\ref{rfe}) has, as a polynomial in $\lambda$, the roots $\lambda=k-\frac{n}{2}$
for $k=1,\dots ,i$ and its leading coefficient is zero. 
Let us first consider $\lambda=i-\frac{n}{2}$, so get after substitution
\begin{eqnarray}
& & A_{i,j}(i-\frac{n}{2},\mu)=
\frac{(-1)^j(i+j-N-1)\dots (-N)}{2^{i+j}i!j!}\cdot
 \\  \nonumber
& & 
\frac{(j-N-\frac{n}{2})\dots (i-N-\frac{n}{2}+1)(1-\frac{n}{2}-N)_i(\mu-N+1)_i}
{(j-\frac{n}{2}-\mu)\dots (1-\frac{n}{2}-\mu)(\lambda+\frac{n}{2}-1)\dots (\lambda+\frac{n}{2}-i+1)}\bigl\lvert_{\lambda=i-\frac{n}{2}}
\end{eqnarray}
and hence
\begin{eqnarray}
A_{i,j+1}(i-\frac{n}{2},\mu)=A_{i,j}(i-\frac{n}{2},\mu)\cdot
\frac{(-1)(i+j-N)(j-N-\frac{n}{2}+1)}{2(j+1)(j-\frac{n}{2}-\mu+1)}.
\end{eqnarray}
Taken together, there remain just two contributions on the left hand side of (\ref{rfe}) given by
$A_{i,j}(i-\frac{n}{2},\mu)$, $A_{i,j+1}(i-\frac{n}{2},\mu)$. Up to a common rational
factor, their sum is proportional to
\begin{eqnarray}
& & i(\frac{n}{2}-1)+(j-N)(-j+N-i+\frac{n}{2}-1)+
\nonumber \\ \nonumber
& & (j+1)(n+2\mu-2-2j)(-1)\frac{(i+j-N)(j-N-\frac{n}{2}+1)}{2(j+1)(j-\frac{n}{2}-\mu+1)}=0,
\end{eqnarray}
which proves the claim. The proof of triviality of the left hand side at special values  
$\lambda=i-1-\frac{n}{2},\dots ,1-\frac{n}{2}$ is completely analogous.

 Note that there are some other equally convenient choices for $\lambda,\mu$ allowing 
the triviality check for (\ref{rfe}), for example based on the choice $\lambda=k+N-1$, $k=1,\dots ,i$
or $\mu=N-k$, $k=1,\dots ,i$.  

The remaining task is to find the leading coefficient on the left hand side of (\ref{rfe})
as a polynomial in $\lambda$. Because 
\begin{eqnarray} 
& & 
(\lambda+\frac{n}{2}-i)_{i-k}\stackrel{\lambda\to\infty}{\sim}\lambda^{i-k},
\nonumber \\
& & 
(\lambda-N+1-k)_{k}\stackrel{\lambda\to\infty}{\sim}\lambda^k,
\end{eqnarray}
the polynomial is of degree $\lambda^{j-i}\frac{\lambda^i}{\lambda^i}=\lambda^{j-i}$, 
$j\geq i$. The leading coefficient of $A_{i,j}(\lambda ,\mu)$ is
\begin{eqnarray} 
\lim _{\lambda\to \infty}\frac{A_{i,j}(\lambda ,\mu)}{\lambda^{j-i}} &=&
(\sum_{k=0}^i(-1)^k {i\choose k}(j-i+1+k)_{i-k}(\mu-N+1)_k)\cdot
\nonumber \\
& & \frac{(-1)^{i+j}(i+j-N-1)\dots (-N)}{2^{i+j}i!j!\,\,(j-\frac{n}{2}-\mu)\dots (1-\frac{n}{2}-\mu)}.
\end{eqnarray}
There are three contributions to (\ref{rfe}):
\begin{eqnarray} \label{limitlambda}
& & (N-j+i)\lim _{\lambda\to \infty}\frac{A_{i,j}(\lambda ,\mu)}{\lambda^{j-i}},
\nonumber \\
& & (j+1)(n+2\mu-2-2j)\lim _{\lambda\to \infty}\frac{A_{i,j+1}(\lambda ,\mu)}{\lambda^{j+1-i}},
\nonumber \\
& & \frac{1}{2}(i+j-N-1)(-i+3j+N-2\mu)\lim _{\lambda\to \infty}\frac{A_{i-1,j}(\lambda ,\mu)}{\lambda^{j+1-i}},
\end{eqnarray}
whose sum is a polynomial in $\mu$ multiplied by common product of linear polynomial. In 
order to prove triviality of this polynomial, it suffices as in the first part of the proof to find 
sufficient amount of its roots and to prove the triviality of its leading coefficient. For 
example in the case $\mu=N-1$, we get from (\ref{limitlambda}) that the coefficients of this 
polynomial are proportional to the sum 
$$
(N-j+i)+\frac{(j+1)(i+j-N)}{(j-i+1)}-\frac{i(-i+3j+N-2(N-1))}{(j-i+1)},
$$  
which equals to zero. The verification of the required property 
for $\mu=N-k$, $k=2,\dots ,i$ is completely analogous. 
This completes the proof.

\hfill
$\square$

Based on the notation \eqref{eqn:sol}, this completes the description of 
\begin{eqnarray}
\mbox{\rm Sol}({\gog}\oplus\gog,\mathrm{diag}({\gog});\mC_{\lambda,\mu}),
\end{eqnarray}
the space of scalar valued singular vectors in the Fourier image of generalized Verma modules 
characterizing solution space of a diagonal branching problem for $so(n+1,1,\mR)$.

It is an interesting observation that the four term functional equation (\ref{rfe})
for $A_{i,j}(\lambda,\mu)$ can be simplified
using the generalized hypergeometric function ${}_3F_2$, 
$$
{}_3F_2(a_1,a_2,a_3;b_1,b_2;z):=
\sum_{m=0}^\infty\frac{(a_1)_{m}(a_2)_{m}(a_3)_{m}}{(b_1)_{m}(b_2)_{m}}\frac{z^m}{m!},
$$
where $a_1,a_2,a_3\in\mC,\, b_1,b_2\in\mC\setminus\{-\mN\}$ and $(x)_m=x(x+1)\dots (x+m-1)$. 
In particular, it can be converted into the four term functional equation 
\begin{eqnarray}
& & \frac{(n+2\lambda)\Gamma(i+j-N)\Gamma(-\frac{n}{2}-\lambda)\Gamma(1-i+j-N+\lambda)\Gamma(1-\frac{n}{2}-\mu)}
{2^{i+j}(-1)^{i+j}\Gamma(1+i)\Gamma(-N)\Gamma(1-N+\lambda)\Gamma(1+j-\frac{n}{2}-\mu)}\cdot
\nonumber \\
& &  
(i(-2j+n+2\mu){}_3F_2(1-i,N-\lambda,1-N+\mu;1-i+j,1-\frac{n}{2}-\lambda;1)+
\nonumber \\
& & 
i(-1+i-j+N-\lambda)(i-3j-N+2\mu)\cdot
\nonumber \\
& & \cdot{}_3F_2(1-i,N-\lambda,1-N+\mu;2-i+j,1-\frac{n}{2}-\lambda;1)+
\nonumber \\
& & 
(i^2-i(-1+n+\lambda)+(j-N)(1+j-N+\lambda))\cdot
\nonumber \\
& & \cdot{}_3F_2(-i,N-\lambda,1-N+\mu;1-i+j,1-\frac{n}{2}-\lambda;1)+
\nonumber \\
& & 
(1+j)(i+j-N)(-1+i-j+N-\lambda)\cdot 
\nonumber \\
& & \cdot{}_3F_2(-i,N-\lambda,1-N+\mu;2-i+j,1-\frac{n}{2}-\lambda;1))=0.
\end{eqnarray}
This functional equation seems not to be present in any standard 
textbook on special function theory of several variables, cf. 
\cite{erd}, \cite{bai}. However, it is a consequence of 
\begin{lemma} 
The generalized hypergeometric function ${}_3F_2$ fulfills
\begin{eqnarray}
& & (j-i+1)_{i}(-i+\lambda+\frac{n}{2})_{i}\,{}_3F_2(-i, N-\lambda,1+\mu-N;1-i+j,1-\frac{n}{2}-\lambda;1)
\nonumber \\
& & = \sum_{k=0}^i(-1)^k {i\choose k}(j-i+1+k)_{i-k}(\lambda+\frac{n}{2}-i)_{i-k}(\mu-N+1)_k(\lambda-N+1-k)_k .
\nonumber \\
\end{eqnarray}
In particular, the diagonal coefficients $A_{i,i}(\lambda,\mu)$ can be written as  
\begin{eqnarray}\label{aii}
& & A_{i,i}(\lambda,\mu)=\frac{\Gamma(2i-N)}{2^{2i+1}i!\, 
\Gamma(-N)\Gamma(1+i-\frac{n}{2}-\mu)\Gamma(1+i-\frac{n}{2}-\lambda)}\cdot
\nonumber \\
& & (\Gamma(1+i-\frac{n}{2}-\lambda)\Gamma(1-\frac{n}{2}-\mu)
{}_3F_2(-i,N-\lambda,1-N+\mu;1,1-\frac{n}{2}-\lambda;1)+
\nonumber \\
& & \Gamma(1+i-\frac{n}{2}-\mu)\Gamma(1-\frac{n}{2}-\lambda)
{}_3F_2(-i,1-N+\lambda,N-\mu;1,1-\frac{n}{2}-\mu;1))
\nonumber 
\end{eqnarray}
\end{lemma}
{\bf Proof:}

The first claim is equivalent to Saalschutz's theorems, cf. \cite{erd}, \cite{bai}. 

As for the second claim, it follows from the definition of ${}_3F_2$ that 
\begin{eqnarray}
& & 
\frac{1}{\Gamma(1+i-\frac{n}{2}-\mu)\Gamma(1+i-\frac{n}{2}-\lambda)}\cdot
\nonumber \\
& & (\Gamma(1+i-\frac{n}{2}-\lambda)\Gamma(1-\frac{n}{2}-\mu)
{}_3F_2(-i,N-\lambda,1-N+\mu;1,1-\frac{n}{2}-\lambda;1)+
\nonumber \\
& & \Gamma(1+i-\frac{n}{2}-\mu)\Gamma(1-\frac{n}{2}-\lambda)
{}_3F_2(-i,1-N+\lambda,N-\mu;1,1-\frac{n}{2}-\mu;1))=
\nonumber \\
& & \sum_{m=0}^i(
\frac{(-i)_m(N-\lambda)_m(1-N+\mu)_m}{(1)_m(1-\frac{n}{2}-\lambda)_m(1-\frac{n}{2}-\mu)_{i-1}}
\nonumber \\
& & +\frac{(-i)_m(N-\mu)_m(1-N+\lambda)_m}{(1)_m(1-\frac{n}{2}-\mu)_m(1-\frac{n}{2}-\lambda)_{i-1}})\frac{1}{m!}.
\end{eqnarray}
Using basic properties of the Pochhammer symbol, e.g. $(x)_m=(-1)^m(-x+m-1)_m$, an 
elementary manipulation yields the result.

\hfill
$\square$

Let us mention that the diagonal coefficients 
$A_{i,i}(\lambda,\mu)=A_{i,i}(\mu,\lambda)$ are, up to a rational multiple coming from the 
ratio of the product of $\Gamma$-functions, symmetric with respect to 
$\lambda\longleftrightarrow \mu$. As a consequence, these polynomials belong to 
the algebra of $\mZ_2$-invariants: 
$$
\mC[\lambda,\mu]^{\mZ_2}\stackrel{\sim}{\to}\mC[\lambda\mu,\lambda+\mu].
$$  
\begin{example}
As an example, in the case of $i=1$ we have 
\begin{eqnarray}
A_{1,1}(\lambda,\mu)=\frac{N(N-1)(\lambda\mu-N(\lambda+\mu)+(1-\frac{n}{2}+N(N-1)))}{(2\lambda+n-2)(2\mu+n-2)},
\end{eqnarray}
and
\begin{eqnarray}
A_{1,j}(\lambda,\mu)&=&
\frac{\Gamma(1+j-N)\Gamma(j-N+\lambda)\Gamma(1-\frac{n}{2}-\mu)}
{2^{j+1}(-1)^{j+1}\Gamma(1+j)\Gamma(-N)\Gamma(1-N+\lambda)\Gamma(1+j-\frac{n}{2}-\mu)}\cdot
\nonumber \\
& & \cdot\frac{(j(-2+n+2\lambda)+2(N-\lambda)(1-N+\mu))}{(n+2\lambda-2)}
\end{eqnarray}
for all $j\in\{1,\dots ,N\}$.
\end{example}
Let us also remark that for special values $\lambda,\mu\in \{m-\frac{n}{2}\,|\, m\in\mN\}$, 
the formula $A_{i,j}(\lambda,\mu)$ simplifies due to the factorization of the 
underlying polynomial. This factorization indicates so called factorization identity, when a 
homomorphism of generalized Verma modules quotients through a homomorphism of 
generalized Verma modules of one of its summands (in the source) or a target homomorphism 
of generalized Verma modules. This naturally leads to the question of full composition 
structure of the branching problem, which goes beyond the formulation in terms of the 
Grothendieck group of the Bernstein-Gelfand-Gelfand parabolic category ${\fam2 O}^\gop$.  

Let us summarize our results.
\begin{theorem}\label{algconstrsingvect}
Let $\gog_{\mathbb R}=so(n+1,1,\mR)$ and $\gop_{\mathbb R}$ its conformal parabolic 
subalgebra with commutative nilradical. Then the diagonal branching problem for 
scalar generalized Verma ${\fam2 U}(\gog\oplus\gog)$-modules induced from characters 
$\chi_{\lambda ,\mu}$ is determined, in the 
Grothendieck group $K({\fam2 O}^{\gop})$ of Bernstein-Gelfand-Gelfand parabolic category 
${\fam2 O}^{\gop}$, by ${\fam2 U}(\gog)$-isomorphism  
in Corollary \ref{charform}, equation (\ref{diag-branch-tk}).

Assuming that $\lambda,\mu\in\mC\setminus\{m-\frac{n}{2}\, |\, m\in\mN\}$,
the summand ${\fam2 M}_{\lambda+\mu-2N}(\gog,\gop)$ in (\ref{diag-branch-tk}) 
is generated by scalar valued singular vector of homogeneity $2N$ of the form 
\eqref{non-homsol}:
$$
p(r,s,t)=\sum_{0\leq i,j,k\leq N|i+j+k=N}A_{i,j}(\lambda,\mu)s^it^jr^k,
$$
where the coefficients $A_{i,j}(\lambda,\mu)$ are given by equation (\ref{solut}). 

In particular, these singular vectors are non-zero, linearly independent 
and of expected weight (induced by the homogeneity), and 
the cardinality of the set of singular vectors is as
predicted by Corollary \ref{charform}. 
\end{theorem}

We also remark that our results for the coefficients $A_{i,j}(\lambda,\mu)$ of 
conformally invariant bilinear differential operators can be directly compared 
to the coefficients $c_{r,s,t}$ derived in \cite{OvRe}. For example, the substitution 
for $\lambda$ and $\mu$, respectively, into our formulas the expression $-n\lambda$
and $-n\mu$, respectively, identifies our linear and quadratic solutions 
in Examples \ref{linear-example} and \ref{quadratic-example} with those given in 
\cite{OvRe} up to the multiple $-1$. In fact, after a tedious but
straightforward computation there is analogous comparison result for all coefficients,
cf. $\cite{OvRe},\, \mbox{page}\, 26, (4.4)$. The reason for different normalizations comes 
from exploiting different initial approaches to the same problem.

\section{Application - the classification of bilinear conformally equivariant differential operators
on line bundles}

Let $M$ be a smooth (complex) manifold equipped with a filtration of its tangent bundle 
$$
0\subset T^1M\subset\dots\subset T^{m_0}M=TM ,
$$
${\fam2 V}\to M$ a smooth (holomorphic) vector bundle on $M$ 
and $J^k{\fam2 V}\to M$ the weighted jet bundle over $M$ defined by 
\begin{eqnarray}
J^k{\fam2 V}=\bigcup_{x\in M} J_x^k{\fam2 V},\, J_x^k{\fam2 V}\stackrel{\sim}{\to} 
\oplus_{l=1}^k\Hom({\fam2 U}_l(gr(T_xM)),\Gamma ({\fam2 V}_x)),
\end{eqnarray}
where ${\fam2  U}_l(gr(T_xM))$ is the subspace of homogeneity at most $l$
elements in 
the universal enveloping algebra of the associated graded algebra $gr(T_xM)$.
A bilinear differential pairing between sections of the bundle ${\fam2 V}$ and 
sections of the bundle ${\fam2 W}$ to sections of the bundle ${\fam2 Y}$ is a 
vector bundle homomorphism 
\begin{eqnarray}
B: J^k{\fam2 V}\times J^l{\fam2 W}\to {\fam2 Y}.
\end{eqnarray}

In the case when $M = G_{\mathbb R}/P_{\mathbb R}$ is a generalized flag manifold, a differential 
pairing is called equivariant if it commutes with the action of $G_{\mathbb R}$ on sections 
of homogeneous vector bundles ${\fam2 V}, {\fam2 W}, {\fam2 Y}$. Denoting 
$\mV$, $\mW$, $\mY$ the inducing complex $P_{\mathbb R}$-representations of  
${\fam2 V}, {\fam2 W}, {\fam2 Y}$,
the space of $G_{\mathbb R}$-equivariant differential pairings is in bijection with
\begin{eqnarray}\label{duality}
& & (({\fam2 U}(\gog)\otimes{\fam2 U}(\gog))\otimes_{{\fam2 U}(\gop)\otimes{\fam2 U}(\gop)}
\Hom(\mV\otimes\mW,\mY))^{P_{\mathbb R}}\, \simeq
\nonumber \\
& & \Hom_{{\fam2 U}(\gog)}({\fam2 M}(\gog,\gop,\mY^\ch),
{\fam2 M}(\gog\oplus\gog,\gop\oplus\gop,(\mV^\ch\otimes\mW^\ch)),
\end{eqnarray}
where the superscript denotes the space of $P_{\mathbb R}$-invariant elements and 
$\mV^\ch, \mW^\ch$ denote the complex dual representations.
\begin{theorem}\label{bilconfdensclass}
Let $G_{\mathbb R}=SO_o(n+1,1,\mR)$ and $P_{\mathbb R}$ its conformal parabolic subgroup, 
$\lambda,\mu\in\mC\setminus\{m-\frac{n}{2}\, |\, m\in\mN\}$ and $N\in\mN$. 
Let us denote by ${\fam2 L}_{\lambda}$ the 
homogeneous line bundle on $n$-dimensional conformal sphere $G_{\mathbb R}/P_{\mathbb R}\simeq S^n$ induced from 
the complex character 
$\chi_{\lambda}$ of $P_{\mathbb R}$. We denote by $\iota :G_{\mathbb R}/P_{\mathbb R}\hookrightarrow G_{\mathbb R}/P_{\mathbb R} 
\times G_{\mathbb R}/P_{\mathbb R}$ the diagonal 
embedding, and by $\iota^\star$ the induced pull-back of sections of vector bundles given by 
restriction on the diagonal.  
Then there exists a set of bilinear conformally equivariant operators 
\begin{eqnarray}
B_N: C^\infty( G_{\mathbb R}/ P_{\mathbb R},{\fam2 L}_{\lambda})\times C^\infty(G_{\mathbb R}/P_{\mathbb R},{\fam2 L}_{\mu})\to
C^\infty(G_{\mathbb R}/P_{\mathbb R},{\fam2 L}_{\lambda+\mu-2N})
\end{eqnarray}
of the form    
\begin{eqnarray}
B_N=\sum_{0\leq i,j,k\leq N|i+j+k=N}A_{i,j}(-\lambda ,-\mu)\iota^\star\, \tilde{s}^i\tilde{t}^j\tilde{r}^k,
\end{eqnarray}
where the coefficients $A_{i,j}(\lambda ,\mu)$ are given by (\ref{solut}) and
\begin{eqnarray}
\tilde{s}=\sum_{i=1}^n\partial^2_{x_i}=\triangle_x,\,\,
\tilde{t}=\sum_{i=1}^n\partial^2_{y_i}=\triangle_y,\,\,
\tilde{r}=\sum_{i=1}^n\partial_{x_i}\partial_{y_i}.
\end{eqnarray}
The set $\{B_N\}_{N\in\mN}$ determines uniquely the set of 
all scalar valued conformally equivariant bilinear differential operators.
\end{theorem}
{\bf Proof:} 

The proof is a direct consequence of Theorem \ref{algconstrsingvect} and duality 
(\ref{duality}), together with the application of inverse Fourier 
transform 
$$
x_j\longleftrightarrow -{\mathbf i}\partial_{\xi_j},\, \partial_{x_j}\longleftrightarrow -{\mathbf i}\xi_j
$$
with ${\mathbf i}\in\mC$ the imaginary unit.

\hfill
$\square$

In many applications, it is perhaps more convenient to express the bilinear differential operators 
in terms of tangent and normal coordinates $t_i=\frac{1}{2}(\xi_i+\nu_i)$ resp. 
$n_i=\frac{1}{2}(\xi_i-\nu_i)$, $i=1,\dots ,n$ to the diagonal submanifold $\iota(G_{\mathbb R}/P_{\mathbb R})\subset 
G_{\mathbb R}/P_{\mathbb R}\times G_{\mathbb R}/P_{\mathbb R}$, where
\begin{eqnarray}
& & r=\frac{1}{4}(\sum_{i=1}^n t_i^2 -\sum_{i=1}^n n_i^2),
\nonumber \\
& & s=\frac{1}{4}(\sum_{i=1}^n t_i^2 +\sum_{i=1}^n n_i^2 +2\sum_{i=1}^n t_in_i),
\nonumber \\
& & t=\frac{1}{4}(\sum_{i=1}^n t_i^2 +\sum_{i=1}^n n_i^2 -2\sum_{i=1}^n t_in_i).
\end{eqnarray} 

\vspace{0.1cm}

\flushleft{{\em Acknowledgment}: 
The present article is a result of a long-lasting cooperation with T. Kobayashi, B. {\O}rsted and
V. Sou\v{c}ek. The author gratefully acknowledges support by the Grant Agency of Czech Republic through 
the grant GACR 22-00091S. 



\vspace{0.3cm}

Petr Somberg

Mathematical Institute of Charles University,

Sokolovsk\'a 83, Praha 8 - Karl\'{\i}n, Czech Republic, 

E-mail: somberg@karlin.mff.cuni.cz.

\end{document}